\documentclass{article}
\newtheorem{thm}{Theorem}[section]
\newtheorem{cor}[thm]{Corollary}

\newtheorem{Def}[thm]{Definition}
\newtheorem{rem}[thm]{Remark}

\newtheorem{ex}{Example}[section]

\newcommand{\be}{\begin{equation}}
\newcommand{\ee}{\end{equation}}
\newcommand{\ben}{\begin{enumerate}}
\newcommand{\een}{\end{enumerate}}

\newcommand{\qed}{\hspace*{\fill}Q.E.D.}  
\title{\large \bf A symmetric Finsler space with Chern connection}
\author{Dariush Latifi,Asadollah Razavi}

\date{}

\begin{document}
\maketitle

\begin{abstract}
We define a symmetry  for a Finsler space with Chern connection
and investigate its implementation and properties  and  find a
relation between them and flag curvature.
\end{abstract}
\section{Introduction}
It is well known that a Riemannian space is locally symmetric if
and only if $\nabla R=0$, for Levi-Civita connection.Therefore we
define a symmetric Finsler space to be a Finsler space whose
hh-curvature is parallel with respect to the Chern connection. As
 the hh-curvature is the Riemannian curvature in Riemannian case,
 this generalizes the definition of Riemannian symmetric space to
Finsler spaces.

\section{Preliminaries}

We will follow verbatim the notation of [1].Let $M$ be a manifold
of dimension $n$, a local system of coordinate $(x^{i}),i=1...n$
on $M$ gives rise to a local system of coordinate $(x^{i},y^{i})$
on the tangent bundle $TM$ through
$y=y^{i}\frac{\partial}{\partial x^{i}}$.
\\

A Finsler structure on $M$ is a function $F:TM\rightarrow
[0,\infty)$ satisfying the following condition :
\begin{description}
    \item[(i)] $F$ is differentiable away from the origin.
    \item[(ii)] $F$ is homogeneous of degree one in $y$ i.e for
    all $\lambda > 0$ $$F(x,\lambda y)=\lambda F(x,y)$$
    \item[(iii)] the $n\times n$ matrix $$g_{ij}=\frac{1}{2}\frac{\partial^{2}F^{2}}{\partial y^{i}\partial
    y^{j}}$$ is positive-definite at every point of $ TM\backslash0
    $.
\end{description}

Denote  the natural projection $TM\backslash 0\rightarrow M$ by
$\pi$. The pullback bundle of $T^{*}M$ is defined by the
commutative diagram
$$\pi^{\ast}(T^{\ast}M)\longrightarrow
  T^{\ast}M  $$
   $$\downarrow\hskip2.3cm\downarrow $$ $$ TM\backslash 0\ \hskip.5cm \longrightarrow \hskip-.5cm^{\pi}\hskip.5cmM $$

The components $g_{ij}$ in $(iii)$ define a section
$g=g_{ij}dx^{i}\otimes dx^{j}$ of the pulled back bundle
$\pi^{\ast}(T^{\ast}M)\otimes\pi^{\ast}(T^{\ast}M)$, where $g$ is
called the fundamental tensor, and usually depends on both $x,y$.

To simplify the computation, we introduce adapted bases for the
bundles $T^{*}(TM\setminus0)$ and $T(TM\setminus0)$, these are :
\\

$\{dx^{i},\frac{\delta
y^{i}}{F}=\frac{1}{F}(dy^{i}+N^{i}_{m}dx^{m})\}$
\hskip1cm,\hskip1cm $\{\frac{\delta}{\delta
x^{i}}=\frac{\partial}{\partial
x^{i}}-N^{m}_{i}\frac{\partial}{\partial
y^{m}},F\frac{\partial}{\partial y^{i}}\}$ \\\\
where
$$N^{i}_{m}=\frac{1}{4}\frac{\partial}{\partial y^{m}}(g^{is}(\frac{\partial g_{sk}}{\partial x^{j}}+\frac{\partial g_{sj}}{\partial x^{k}}-\frac{\partial g_{kj}}{\partial x^{s}})y^{j}y^{k})$$
In fact they are dual to each other. The vector space spanned by
$\frac{\delta}{\delta x^{i}}( resp.\  F\frac{\partial}{\partial
y^{i}})$ is called \emph{horizontal(resp.\ vertical)} subspace of
$T(TM)$.

Let $(M,F)$ be a Finsler manifold. There exists a symmetric
connection $$\nabla:\Gamma(T(TM))\times
      \Gamma(\pi^{\ast}(TM)) \longrightarrow
      \Gamma(\pi^{\ast}(TM))$$
whose Christoffel symbols are given by : $$\Gamma^{i}_{jk} =
\frac{g^{is}}{2}(\frac{\delta g_{sj}}{\delta x^{k}}-\frac{\delta
g_{jk}}{\delta x^{s}}+\frac{\delta g_{ks}}{\delta x^{j}})$$this
connection is called Chern connection and has the following
properties :
\begin{description}
    \item[(i)] the connection 1-form does not depend on $dy$;
    \item[(ii)] the connection $\nabla$ is almost $g-$compatible
    in the sense that;\\
    $$\nabla_{\delta / \delta x^{s}}^{\hskip.35cm
      g_{ij}}=0 \hskip1.5cm and \hskip1.5cm \nabla_{F \frac{\partial}{\partial y^{s}}}^{\hskip.35cm
      g_{ij}}=2A_{ijs}$$
\end{description}
where $A_{ijs}$ is the component of Cartan tensor.
\\

The curvature of Chern connection can be splitted into two
components according to the vector argument being horizontal or
vertical. The first is the hh-curvature tensor \\$$(\nabla_{\delta
/ \delta x^{k}} \nabla_{\delta / \delta x^{l}}-\nabla_{\delta /
\delta x^{l}} \nabla_{\delta / \delta x^{k}}-\nabla_{[\delta /
\delta x^{l} , \delta / \delta x^{k}] })\partial/\partial
x^{j}=R^{i}_{jkl} \partial/\partial x^{i}$$ \\ where

$$ R^{i}_{jkl}=\frac{\delta \Gamma^{i}_{jl} }{\delta
   x^{k}}-\frac{\delta \Gamma^{i}_{jk} }{\delta x^{l}}+\Gamma^{i}_{hk}\Gamma^{h}_{jl}-\Gamma^{i}_{hl}\Gamma^{h}_{jk}$$
\\
the second is hv-curvature tensor
\\
$$(\nabla_{F\frac{\partial}{\partial
   y^{k}}}\nabla_{\delta /\delta x^{l}}-\nabla_{\delta /\delta
   x^{l}}\nabla_{F\frac{\partial}{\partial y^{k}}}-\nabla_{[F\frac{\partial}{\partial y^{k}},\frac{\delta}{\delta
   x^{l}}]})\partial/\partial x^{j}=P^{i}_{jkl}\partial/\partial
   x^{i}
   $$
\\
where $$ P^{i}_{jkl}=-F\frac{\partial \Gamma^{i}_{jk}}{\partial
   y^{l}}$$

\section{Symmetric Finsler Space}
\begin{Def}
A Finsler space throughout which the hh-curvature $R$ possesses
vanishing covariant derivative with respect to horizontal vector
field will be called a symmetric Finsler space :
$$\nabla_{\delta /\delta x^{l}}R^{i}_{jkh}=0$$
\end{Def}
\begin{ex}
Evidently a locally Minkowski space is symmetric space, as the
Chern connection coefficients $\Gamma^{i}_{jk}$ vanish
identically.
\end{ex}

Let
$$R^{i}_{kl}=l^{j}R^{i}_{jkl}\hskip1cm and \hskip1cm R^{i}_{j}=l^{k}R^{i}_{}jk$$
where $l^{i}=\frac{y^{i}}{F}$, then we have :
\begin{equation}
  \nabla_{h}R^{i}_{kl}=0
\end{equation}
and  $$\nabla_{h}R^{i}_{j}=0$$ where the index $h$ denotes
$\delta/\delta x^{h}$.
\\

Moreover   the horizontal covariant derivative of  the following
contraction will be zero.$$R^{i}_{ikh}=R_{kh}\hskip.5cm,\hskip.5cm
R^{i}_{ik}=R_{k}\hskip.5cm,\hskip.5cmR^{i}_{i}=(n-1)R$$ In order
to obtain further consequences in a symmetric Finsler space we
compute $$R^{i}_{hkl}=\frac{\partial}{\partial
y^{h}}R^{i}_{kl}+y^{j}\frac{\partial}{\partial
y^{h}}[\dot{A}^{i}_{jl|k}-\dot{A}^{i}_{jk|l}+\dot{A}^{s}_{jl}\dot{A}^{i}_{sk}-\dot{A}^{s}_{jk}\dot{A}^{i}_{sl}]$$
where $\dot{A}^{i}_{jl|k}$ is the horizontal covariant derivative
of $\dot{A}^{i}_{jl}$ with respect to $\frac{\delta}{\delta
x^{k}}$ and $\dot{A}^{i}_{jl}=A^{i}_{jl|s}l^{s}$ and denote
$$ y^{j}\frac{\partial}{\partial
y^{h}}[\dot{A}^{i}_{jl|k}-\dot{A}^{i}_{jk|l}+\dot{A}^{s}_{jl}\dot{A}^{i}_{sk}-\dot{A}^{s}_{jk}\dot{A}^{i}_{sl}]$$
by $D^{i}_{hkl}$, thus we
have$$R^{i}_{hkl}=\frac{\partial}{\partial
y^{h}}R^{i}_{kl}+D^{i}_{hkl}$$
then$$\nabla_{p}R^{i}_{hkl}=\nabla_{p}\frac{\partial}{\partial
y^{h}}R^{i}_{kl}+\nabla_{p}D^{i}_{hkl}$$and if
\begin{equation}
\nabla_{p}D^{i}_{hkl}=0
\end{equation}
then$$\nabla_{p}R^{i}_{hkl}=\nabla_{p}\frac{\partial}{\partial
y^{h}}R^{i}_{kl}$$ Now we can prove:
\begin{thm}
Let $M$ be a Finsler space for which $\nabla_{p}D^{i}_{hkl}=0$.
Then $M$ is symmetric if and only if
\begin{description}
    \item[(i)] $\nabla _{h}R^{i}_{kl}=0$
    \item[(ii)] $\frac{\partial R^{i}_{kl}}{\partial
y^{m}}\dot{A}^{m}_{hp}+R^{m}_{kl}\Gamma^{i}_{mph}-R^{i}_{ml}\Gamma^{m}_{kph}-R^{i}_{km}\Gamma^{m}_{lph}=0$
\end{description}
\end{thm}
\emph{Proof}: For a symmetric space we have $\nabla
_{h}R^{i}_{kl}=0$ and if moreover $\nabla_{p}D^{i}_{hkl}=0$ holds
then with the help of the commutation formula :
\begin{equation}
\frac{\partial}{\partial
y^{h}}\nabla_{p}R^{i}_{kl}-\nabla_{p}\frac{\partial}{\partial
y^{h}}R^{i}_{kl}=\frac{\partial R^{i}_{kl}}{\partial
y^{m}}\dot{A}^{m}_{hp}+R^{m}_{kl}\Gamma^{i}_{mph}-R^{i}_{ml}\Gamma^{m}_{kph}-R^{i}_{km}\Gamma^{m}_{lph}
\end{equation}
where $\Gamma^{i}_{mph}$ is used to denote
$\frac{\partial}{\partial y^{h}}\Gamma^{i}_{mp}$ , we have
\begin{equation}
\frac{\partial R^{i}_{kl}}{\partial
y^{m}}\dot{A}^{m}_{hp}+R^{m}_{kl}\Gamma^{i}_{mph}-R^{i}_{ml}\Gamma^{m}_{kph}-R^{i}_{km}\Gamma^{m}_{lph}=0
\end{equation}
Conversely let us suppose that (i) and (ii) satisfy then from the
commutation formula we have the following equation

$$\nabla_{p}R^{i}_{hkl}=-\frac{\partial R^{i}_{kl}}{\partial
y^{m}}\dot{A}^{m}_{hp}-R^{m}_{kl}\Gamma^{i}_{mph}+R^{i}_{ml}\Gamma^{m}_{kph}+R^{i}_{km}\Gamma^{m}_{lph}$$
therefore $\nabla _{p}R^{i}_{hkl}=0$ .\qed
\\
\\ As for a Ladsberg space $\dot{A}_{ijk}=0$ we have the following
corollary.

\begin{cor}
A Landsberg space is symmetric if and only if (1) and (4) hold.
\end{cor}
\begin{Def}
A Finsler structure $F$ is said to be of Berwald type if the Chern
connection coefficient $\Gamma^{i}_{jk}$, in natural coordinates,
have no $y$ dependence.
\end{Def}
Landsberg spaces include Berwald type spaces, now we have:

\begin{thm}
Let $(1)$ holds in a Berwald type space then it is symmetric.
\end{thm}

If $X$ is a nowhere zero vector field defined on an open subset
$O$ of a Finsler manifold $(M,F)$, then we may associate to $X$ a
Riemannian metric on $O$. A particularly interesting case of this
construction is when the integral curves of the vector field are
geodesics of the Finsler metric on $M$.
\begin{Def}
Let $M$ be a Finsler manifold and let $v_{m}\in T_{m}M$ be a
nonzero vector. If $P\subset T_{m}M$ is a two-dimensional subspace
containing $v_{m}$, and $X$ is a geodesic vector field on a
neighborhood of $m$ such that $X(m)=v_{m}$ then the sectional
curvature of the Riemannian metric associated to $X$ at the plan
$P$ is called flag curvature of $M$ at the flag $(P,v_{m})$.
\end{Def}

\begin{thm}
Let $(M,F)$ be a connected Finsler manifold with constant flag
curvature of dimension at least 3, if (2) and (4) hold then $M$ is
symmetric.
\end{thm}
\emph{Proof}: Let $(M,F)$ has constant flag curvature $\lambda$,
and its dimension $n$ is at least $3$, then from [1] we have
$$R^{i}_{kl}=\lambda(\delta^{i}_{k}l_{l}-\delta^{i}_{l}l_{k})$$
thus $\nabla_{h}R^{i}_{kl}=0$.\qed

\begin{cor}
A Landsberg space with constant flag curvature and dimension
greater than 2 which satisfies in (4) is symmetric.
\end{cor}

\begin{cor}
Any Berwald type space with constant flag curvature and dimension
greater than 2 is symmetric.
\end{cor}
\begin{rem}
Traditionally Riemannian symmetric spaces were first defined and
studied by E.Cartan [2],by means of symmetries i.e involutive
isometries fixing a point. It is shown that locally it is
equivalent to having parallel curvature tensor field.
Affine(locally and globally) symmetric space are defined by
replacing isometries with affine transformation, and is equivalent
to having a torsionfree connection with parallel curvature tensor
field. O.Loos considered affine symmetric space as a smooth
manifold with an appropriate operation and produced the required
connection. Therefore it seem reasonable to consider a space
together with a nice operation and find a suitable connection .
\end{rem}

\noindent

Department of Mathematics and Computer Science, Amirkabir
University of Technology,P.O.Box 15875-1433.Tehran,Iran

\emph{E-mail address}: dlatifi@aut.ac.ir , arazavi@aut.ac.ir
\end{document}